\documentclass[pre,preprint,showpacs, showkeys, preprintnumbers,amsmath,amssymb,floatfix]{revtex4}
\usepackage{graphicx}
\usepackage{dcolumn}
\usepackage{bm}

\usepackage[normal]{subfigure}
\usepackage{chemarr}
\usepackage[matrix,curve,arrow,tips,frame]{xy}
\graphicspath{{figures/}}
 \DeclareMathOperator{\sgn}{sign}
\begin{document}

\title{Limit cycles in the presence of convection,\\ a travelling wave analysis}

\author{E.~H.~Flach}
\email{flach@indiana.edu}
\affiliation{Complex Systems, Indiana University 	School of Informatics,\\
1900 East Tenth Street, Eigenmann Hall 906, Bloomington, IN 47406, USA}
\altaffiliation[Also at ]{Centre for Mathematical Biology,
Mathematical Institute, 24-29 St Giles', Oxford, OX1 3LB, UK}
\author{S.~Schnell}
\affiliation{Complex Systems, Indiana University 	School of Informatics,\\
1900 East Tenth Street, Eigenmann Hall 906, Bloomington, IN 47406, USA}
\author{J.~Norbury}
\affiliation{Mathematical Institute, 24-29 St Giles', Oxford, OX1 3LB, UK}

\date{\today}

\begin{abstract}
We consider a diffusion model with limit cycle reaction functions.
In an unbounded domain, diffusion spreads pattern outwards from the source.
Convection adds instability to the reaction-diffusion system.
The result of this instability is a readiness to create pattern.
We choose the Lambda-Omega reaction functions for their simple limit cycle.
We carry out a transformation of the dependent variables into polar form.
From this we consider the initiation of pattern to approximate a travelling wave.
We carry out numerical experiments to test our analysis.
These confirm the premise of the analysis, that the initiation can be modelled by a travelling wave.
Furthermore, the analysis produces a good estimate of the numerical results.
Most significantly, we confirm that the pattern consists of two different types.
\end{abstract}

\maketitle

\section{Introduction}
Morphological patterning, such as animal coat markings, may be caused by a chemical field~\cite{Turing1952The-Chemical-Ba,Murray1981A-Pre-pattern-f,Murray1981On-Pattern-Form}.
The Turing model now has strong experimental support~\cite[for a review]{Maini2006The-Turing-Mode},~\cite{Sick2006WNT-and-DKK-Det,Jung1998Local-inhibitor}.
This mechanism has some limitations, and so we continue to investigate variations of the model~\cite{Maini1996Spatial-and-spa}.

We consider the standard reaction-diffusion model with the addition of convection~\cite{Bamforth2001Flow-distribute,Bamforth2000Modelling-flow-,Kuznetsov1997Absolute-and-co}.
The more general form is a system with advection, as in~\cite{Andresen1999Stationary-spac}.
Our motivation is theoretical, to see the effect of convection on the robustness of the pattern formation.
We consider that a weak amount of convection, or a similar effect, may be present \textit{in vivo}.
This may soften the standard requirements for the formation of pattern, increasing the applicability of the model.

There may be a direct biological application, such as the formation of the vertebral precursors (somitogenesis).
The organism growth could produce such a convective effect~\cite{Kaern2001Chemical-waves-}.
However, convection-induced patterning is not generally supposed for somitogenesis~\cite{Schnell2002Models-for-patt}.

Our system does not fit the standard Turing analysis since our reaction functions are already unstable.
However, we have shown previously that an unstable function may produce a Turing pattern~\cite{Flach2007Turing-pattern-}.
Furthermore, we have shown that, in the presence of convection, we can see similar behaviour irrespective of the stability of the fixed point~\cite{Flach2007Turing-pattern-}.

Experimentally, convection is introduced from the boundaries~\cite{Kaern1999Flow-distribute, Miguez2006Robustness-and-}.
The effect of boundary conditions is also likely to be relevant \textit{in vivo}.
The boundary can have a significant effect on the pattern formed~\cite{Bamforth2000Modelling-flow-}.
However, we follow the suggestion of Cross and Hohenberg to examine the system in a boundary-free environment first~\cite{Cross1993Pattern-formati}.

In numerical simulations we make form an initial point disturbance.
Pattern is formed, spreading outwards from its initiation point.
The disturbance is oscillatory and complex.
Our focus here is the way in which the pattern propagates.

The speed of propagation has been given for a two-species system with no convection~\cite{Bricmont1994Stability-of-mo} and one with equal convection~\cite{Bamforth2000Modelling-flow-}.
A related theoretical study considers the relationship between the onset of the instability and the longer-term behaviour of the non-equilibrium state~\cite{Sherratt1998Invading-wave-f}.
Here the simplest reaction-diffusion system with oscillatory kinetics was considered, and complex behaviours were found.

Limit cycles are inherent in these oscillatory models, since an unstable spiral at the steady state must be bounded~\cite{Schnakenberg1979Simple-chemical}.
In this paper, we seek to exploit the limit cycle to aid our analysis.
To this end, we select the Lambda-Omega reaction functions.
These are chosen to produce the simplest limit cycle possible: the unit circle.

The reaction functions chosen in a full model correspond to actual reaction mechanisms.
Examples of these include the Belousov-Zhabotinsky reaction, and Brusselator-type systems such as the chlorine dioxide-iodine-malonic acid reaction.
The Schnakenberg model can be considered an intermediate step towards these functions.
For an initial investigation, we have chosen simple, related functions.
The intention is to extend the analysis to more realistic models.
Secondarily, the results in this simple case may form a useful basis for comparison to more complex systems.

We transform the problem so that the oscillations are removed, or at least reduced.
This brings the onset of the pattern into much sharper relief.
The two aspects -- the pattern, and the initiation of the pattern -- are clearly distinguished.
Thus this onset of pattern is effectively a phase transition.

The transformation is to convert the dependent variables into polar form.
We then consider that the formation of the limit cycle approximates a travelling wave, and so we employ a Fisher solution to the problem~\cite{Wang1988Exact-and-expli}.
The relative sizes of the parameters affects our analysis: we consider various situations.
We carry out numerical experiments to discover which of these estimates are valid.

\section{Limit cycles and the $\lambda$-$\omega$ function}
We suppose that some form of chemical mechanism is the underlying basis for biological pattern formation.
If the chemicals are well-mixed, then the law of mass action valid, and an ODE is an appropriate model: \begin{eqnarray}
    u' &=& f \ , \nonumber\\
    v' &=& g \ . \label{eq:ODE}\end{eqnarray}

\begin{figure}
\centerline{\includegraphics[width=4.2in]{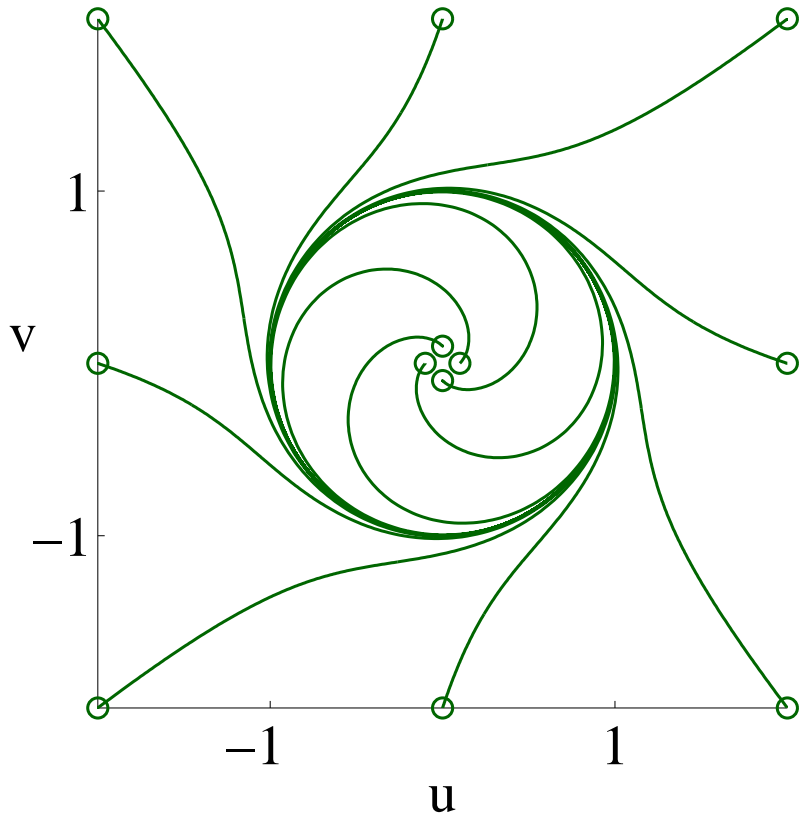}}
\caption{Phase space for the limit-cycle reaction~\eqref{eq:ODE},~\eqref{eq:lambdaomegafunctions}, given by numerical solution.
The phase curves spiral out from the steady state (the origin) to meet the limit cycle (the unit circle).
The trajectories starting far away from the steady state spiral into the limit cycle.
Solid curves starting with circles are the trajectories.} \label{fig:ODE}
\end{figure}
We consider reaction functions chosen for their simplicity in the form of the limit cycle produced:\begin{eqnarray}
    f &=& -v +u(1-u^2-v^2) \ , \nonumber\\
    g &=&  u +v(1-u^2-v^2) \ . \label{eq:lambdaomegafunctions}\end{eqnarray}
This is a simple form of the $\lambda$-$\omega$ class of functions.
Here the steady state is the origin and the limit cycle is the unit circle.
There is no parameter to determine the stability: the steady state is unstable and the limit cycle stable.
This behaviour is clear in the $(u,v)$ phase plane (see~{\sc Figure}~\ref{fig:ODE}).
The circled points on the diagram show the start of the phase plane trajectories.
The phase curves spiral out from the steady state to join the limit cycle.
The trajectories starting far away from the steady state spiral into the limit cycle.

These functions are chemically unrealistic, as they stand.
However, we can relate these functions to the Schnakenberg reactions, chosen to be the simplest chemical form which can produce a limit cycle.
This mechanism was proposed theoretically, but has been used as a model for actual reaction mechanisms~\cite{Epstein1998An-Introduction, Gray1994Chemical-Oscill}.
The key reaction there is the cubic autocatalytic one, $U + 2V \rightarrow 3V$.
Using the law of mass action, this produces the term $uv^2$.
The reverse step gives $v^3$.
These types of terms form the core of the $\lambda$-$\omega$ functions, the cubic terms.

We consider a new coordinate system for the dependent variables
${u, v}$ as a polar form ${r, \theta}$ as follows:\begin{eqnarray}
    r^2 &=& u^2 +v^2            \ , \nonumber\\
    \tan\theta &=&  \frac{v}{u} \ . \end{eqnarray}
By differentiating these identities, we transform the first order ODE~\eqref{eq:ODE} into:\begin{eqnarray}
    \dot{r} &=& r(1-r^2)  \ , \nonumber\\
    \dot{\theta} &=& 1    \ , \label{eq:polaridentities}\end{eqnarray}
The $\theta$ equation clearly decouples and is resolvable.
The remaining single-species ODE in $r$ has steady states at $-1, 0, 1$.
The negative state is unrealistic because $r$ must be positive, zero is unstable and one is stable.
The point $r=1$ is our limit cycle: the unit circle.
Furthermore, $\theta \approx t$, time and so progression around the limit cycle is constant.

This type of analysis is often used in similar cases, such as the analysis of travelling wave trains~\cite{Murray1989Mathematical-Bi}.

\section{Pattern formation}\label{sec:pattern}
Biological pattern formation is by definition spatially differentiated.
There are many models for this.
One is that chemicals diffuse within an organism, then cells respond differently dependent on the  concentration of one of these chemicals.
If the concentration of the chemicals has formed a pattern, then this is reflected by the cells.

For an initial analysis, we choose the simplest case: one spatial dimension and some diffusion.
By diffusion we refer to the averaged gross effect of random motion of the chemicals; with passive movement this is equivalent to Brownian motion.
In the case of general reaction functions, this system is the one proposed by Turing~\cite{Turing1952The-Chemical-Ba}.

We add a convective term to the system.
This can be considered a disturbance to the system, as an investigation of stability.
For a small amount of convection, this could correspond to axial growth~\cite{Kaern2001Chemical-waves-}.
However, for convection to be significant in axial growth we would require extremely slow diffusion of the chemicals to fit the model.

The general system is then as follows: \begin{eqnarray}
    u_t &=&  \varepsilon_1 u_{\xi\xi} -p u_\xi+ f              \ , \nonumber\\
    v_t &=&  \varepsilon_2 v_{\xi\xi} -q v_\xi + g \ . \label{eq:original}\end{eqnarray}
We consider different convection on each species ($p \neq q$).
This effect could be created in a chemical flow reactor, with one of the reactants held fixed in a packed bed~\cite{Kaern2001Chemical-waves-}.
For biological applications, one of the chemicals may be held within a cell, while the other flows freely.
If both chemicals flow freely, it is possible that their movement may be hampered by obstacles such as the extra-cellular matrix.
In this case, a larger molecule may be affected more strongly than a smaller one, and different convection speeds may result for each chemical.
This concept parallels that of different diffusion rates, as in the Turing model.

We consider that $p>0$.
If this is not the case, we make the transformation $(p, q, \xi) \rightarrow (-p, -q, -\xi)$, giving a positive value for $p$ with no change to the form of~\eqref{eq:original}.
We can remove one of the convective terms by a simple change of coordinates $x = \xi - pt$: \begin{eqnarray}
    u_t &=&  \varepsilon_1 u_{xx} + f              \ , \nonumber\\
    v_t &=&  \varepsilon_2 v_{xx} - \gamma v_x + g \ , \label{eq:full}\end{eqnarray}
with $f$ and $g$ as in~\eqref{eq:lambdaomegafunctions} and $\varepsilon_1$, $\varepsilon_2$ and $\gamma$ positive constants.
From the original system~\eqref{eq:original} we have $\gamma = q-p$.
If $\gamma$ turns out to be negative, we can reverse the sign as we did for $p$ and $q$.
Equivalently, we choose the coordinate change $x = pt -\xi$ previously, and set $\gamma = p-q$.
In the case where the convection is the same on both species ($p=q$), then both convection terms are removed ($\gamma=0$).
This recovers the basic Turing model.

We have reaction functions that have limit-cycle behaviour and convection, which is known to drive instability.
The appearance of pattern is then to be expected, although the form might be more difficult to predict.
Satnoianu, Merkin and Scott~\cite{Satnoianu1998Spatio-temporal} studied a similar system to~\eqref{eq:full} previously, with Schnakenberg reaction functions in place of the $\lambda$-$\omega$ ones.
They found that periodic behaviour is emergent in the system over a broad parameter range.

\begin{figure}
\centerline{\includegraphics[width=4.2in]{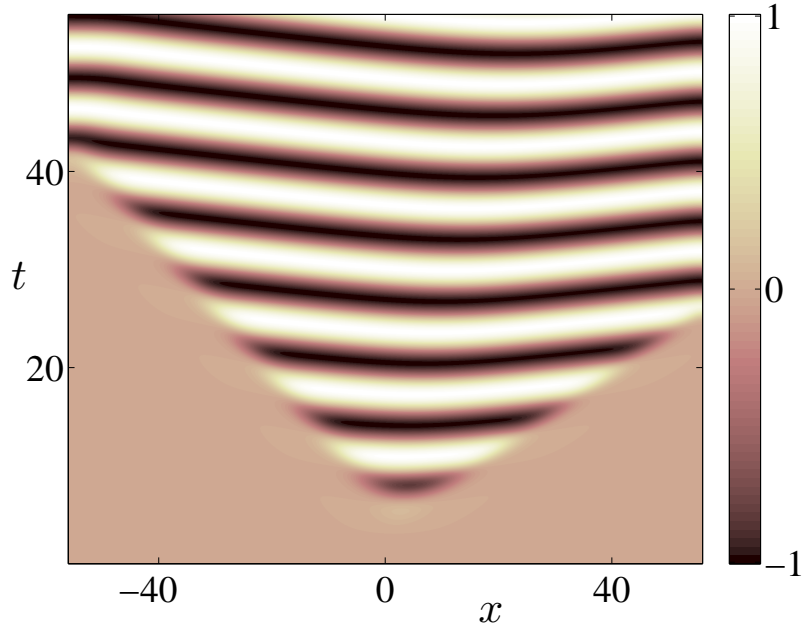}}
\caption{Pattern found for a diffusion system with convection and limit-cycle reaction kinetics~\eqref{eq:full}.
The initial disturbance propagates and becomes pronounced, forming a regular pattern with aligned oscillations.
The propagation is linear, forming a V-shape.
The convective effect is slight: the pattern is skewed slightly to the right.
In this case, the direction of the internal oscillations is distinct from the angle of propagation of the pattern.
This is a numerical solution of~\eqref{eq:full} using NAG D03PCF, plotting species $u$ with $\gamma=\varepsilon_1=\varepsilon_2=1$.
The reactants are initially at steady state: $\left(u,v\right)=\left(0, 0 \right)$, with a small disturbance at $x=0$.
The boundaries are held at zero derivative: $u_x=0$, $v_x=0$.}\label{fig:pattern}
\end{figure}
In the numerical experiment, we start at the steady state $(0, 0)$, except for a small disturbance at $x=0$.
We try to simulate a boundless environment -- to this end we find zero derivative boundary conditions the most effective.
The initial disturbance propagates and becomes pronounced, forming a regular pattern with aligned oscillations.
The propagation is linear, forming a V-shape.
The convective effect is to skew the pattern to the right ({\sc Figure}~\ref{fig:pattern}).
Given the parameters in the figure ($\gamma=\varepsilon_1=\varepsilon_2=1$), the direction of the internal oscillations is distinct from the angle of propagation of the pattern.
The emergence of this pattern is the main focus of our study.

\section{Travelling wave analysis}
Before we start any specific analysis, we rescale the spatial variable $x$ to remove one of the parameters.
We choose $x=\sqrt{\varepsilon_1}y$.
This yields the system: \begin{eqnarray}
    u_t &=&  u_{yy} + f              \ , \nonumber\\
    v_t &=&  v_{yy} + \bar{\varepsilon} v_{yy} - \bar{\gamma} v_y + g \ , \label{eq:rescaled}\end{eqnarray}
where $\bar{\gamma} = \gamma/\sqrt{\varepsilon_1}$.
Following suite, we expect $\bar{\varepsilon_2} = \varepsilon_2/\varepsilon_1$, but we go one step further in defining $\bar{\varepsilon} = \bar{\varepsilon_2} - 1$.
This produces a symmetric first spatial derivative on both equations, with $\bar{\varepsilon}$ the difference between the two diffusion rates.

We wish to again convert ${u, v}$ into the polar form $(r, \theta)$.
During the conversion, we see that there is grouping of the terms parameterised by $\bar{\varepsilon}$ and $\bar{\gamma}$: \begin{eqnarray}
    r_t &=& r_{yy} - r\theta_y^2 + \frac{v}{r}[\bar{\varepsilon}v_{yy} - \bar{\gamma}v_y] + r(1-r^2) \ , \nonumber\\
    \theta_t &=& \theta_{yy} + 2\frac{r_y}{r}\theta_y  + \frac{u}{r^2}[\bar{\varepsilon}v_{yy} - \bar{\gamma}v_y] + 1 \ . \end{eqnarray}
We complete the transformation: \begin{eqnarray}
    r_t &=& \left(1+\bar{\varepsilon}\sin^2\theta\right)r_{yy} \nonumber\\
    && + \left(2\bar{\varepsilon}\sin\theta\cos\theta.\theta_y - \bar{\gamma}\sin^2\theta\right)r_y \nonumber\\
    && + \left(1 -r^2 -\theta_y^2 +\bar{\varepsilon}\sin\theta\cos\theta.\theta_{yy} -\bar{\varepsilon}\sin^2\theta.\theta_y^2 -\bar{\gamma}\sin\theta\cos\theta.\theta_y\right)r \ , \nonumber\\
    \theta_t &=& \left(1 +\bar{\varepsilon}\cos^2\theta\right)\theta_{yy} \nonumber\\
    && +\left(2\frac{r_y}{r} +2\bar{\varepsilon}\frac{r_y}{r}\cos^2\theta -\bar{\varepsilon}\sin\theta\cos\theta.\theta_y -\bar{\gamma}\cos^2\theta \right)\theta_y \nonumber\\
    && +1 +\frac{\bar{\varepsilon}r_{yy} -\bar{\gamma}r_y}{r}\sin\theta\cos\theta \ . \end{eqnarray}

\begin{figure}
\centerline{\includegraphics[width=4.2in]{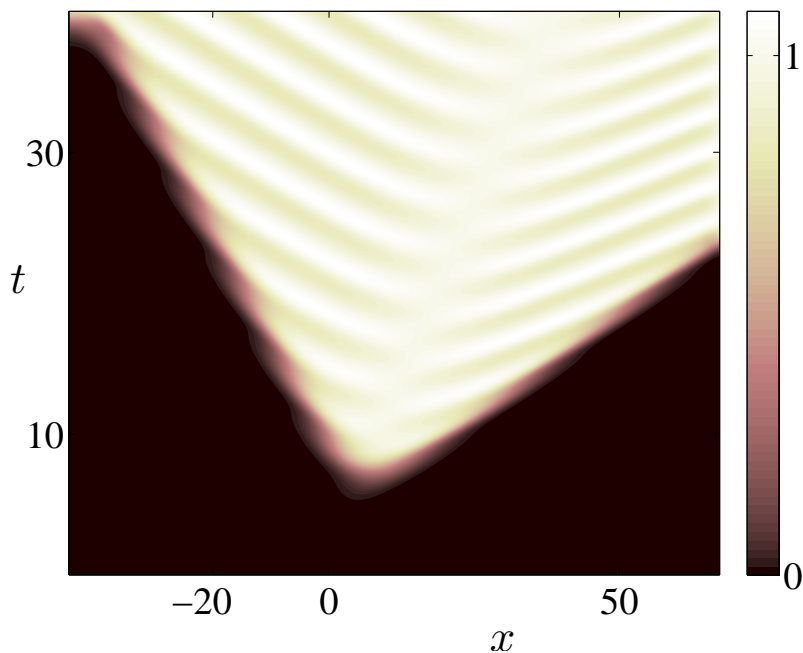}}
\caption{Polar form of the pattern.
The new coordinate $r$ transitions sharply from $0$ to $1$.
This is a travelling wave, propagating outwards.
There are minor ripples in the established solution: the limit cycle settled on by the PDE is not the unit circle, but is close.
The oscillation is at an angle to the travelling wave front.
There is a difference in behaviour between the left and right sides of the pattern: the angle of alignment and the frequency of the pattern is different on either side.
There is also a clear centre to the propagation, roughly at $x= \gamma/2 t$.
This is a numerical solution of~\eqref{eq:full} using NAG D03PCF, plotting species $u$ with $\varepsilon_1=\varepsilon_2=1$, $\gamma=2$.
The reactants are initially at steady state: $\left(u,v\right)=\left(0, 0 \right)$, with a small disturbance at $x=0$.
The boundaries are held at zero derivative: $u_x=0$, $v_x=0$.}\label{fig:polar}
\end{figure}
If we translate our numerical results into this polar form, the onset of the pattern becomes very clearly demarked: the new coordinate $r$ transitions sharply from $0$ to $1$.
The primary behaviour is that of a travelling wave, propagating outwards.
The internal behaviour of the pattern, the steady oscillation, is reduced to a secondary effect.
There are minor ripples in the established solution: the limit cycle settled on by the PDE is not the unit circle, but is close.
There is a difference in behaviour between the left and right sides of the pattern: the angle of alignment and the frequency of the pattern is different on either side ({\sc Figure}~\ref{fig:polar}).
Next we examine this travelling wave analytically.

\subsection{Simple version}
We consider the very simplest situation: equal diffusion on both chemical species ($\bar{\varepsilon}=0$) combined with no convection of $v$ ($\bar{\gamma}$=0).
We make a simplifying assumption for an initial analysis: $\theta_y\approx 0$, $\theta_{yy} \approx 0$.
The problem reduces to \begin{eqnarray}
    r_t &=& r_{yy} + r(1-r^2) \ , \nonumber\\
    \theta_t &=& 1 \ . \end{eqnarray}
Then $\theta$ decouples, as for the ODE.
The solution for $\theta$ is $\theta = t$, to within an arbitrary constant.

The remaining equation in $r$ is of the form similar to the Fisher equation and should therefore yield a propagating wave solution.
This analysis is covered in greater detail in~\cite{Murray1989Mathematical-Bi}.
We look for a solution of the form \begin{eqnarray}
    r(y,t) = R(z) \ , && z=y-\bar{c}t \ . \end{eqnarray}
which gives \begin{eqnarray}
    R'' + \bar{c} R' + R(1-R^2) = 0 \ . \end{eqnarray}
We carry out a phase plane analysis in the $(R, R')$ phase plane.
We find that $\bar{c}>0$ gives a stable point at $(R, R') = (0,0)$, which suggests the solution we are looking for.
The point $R=1$, which corresponds to our limit cycle, is always a saddle point.
We are looking for a phase plane trajectory that leaves this saddle point and goes to the zero
steady state: this will be our travelling wave solution.

Small wave speed, $\bar{c} < 2$, gives a spiral in the phase plane, which would imply $r < 0$ at some point on the trajectory.
From the definition of $r$ we know this to be impossible and thus this is unrealistic for a travelling wave solution.
For $\bar{c} \geq 2$ we have a node and the trajectory discussed above, leaving from the saddle point will head directly to the node, remaining in the fourth quadrant of the phase plane and therefore realistic.
This trajectory equates to the travelling wave that we are looking for.
We expect that the least speed wave will be achieved so we predict a wave of speed $\bar{c}=2$.
Converting this back into our unscaled system, we have $c_R = 2\sqrt{\varepsilon_1}$ to the right.
There is also a solution to the left: $c_L = -2\sqrt{\varepsilon_1}$.
In the original system, we have $\textit{wavespeed} = p \pm 2\sqrt{\varepsilon_1}$.

\subsection{Full system}
We now apply the method for the full system.
Here we may have different diffusion on the two species ($\bar{\varepsilon} \ne 0$) or some convection ($\bar{\gamma} > 0$). 
We again make the assumption $\theta' \approx 0$, $\theta'' \approx 0$.
We introduce the travelling wave coordinate $z= y-ct$ and look for a solution $R(z) = r(y, t)$ .
We also linearise the system, dropping the $R^3$ term:
\begin{eqnarray}
\left(1+\bar{\varepsilon}\sin^2\theta\right)R'' +\left(c  - \bar{\gamma}\sin^2\theta\right)R' + R &=& 0 \ , \nonumber\\
\frac{\bar{\varepsilon}R'' - \bar{\gamma}R'}{R}\sin\theta\cos\theta +1 &=& 0 \ . \end{eqnarray}
Using the same approach as before, we find $\bar{c} = \bar{\gamma}\sin^2\theta \pm 2\sqrt{1+\bar{\varepsilon}\sin^2\theta}$.
Converting this back to the unscaled coordinate system we have $c= \gamma\sin^2\theta \pm 2\sqrt{\varepsilon_1\cos^2\theta+\varepsilon_2\sin^2\theta}$.
We see that $\gamma$ and $\varepsilon_2$ remain linked by the $\sin^2$ term.
These parameters occur on the $v$ differential equation, and $v = r\sin\theta$.

In the original system we have $\textit{wavespeed} = p\cos^2\theta + q\sin^2\theta \pm 2\sqrt{\varepsilon_1\cos^2\theta+\varepsilon_2\sin^2\theta}$.
Again, the parameters for each species $u$ and $v$ are linked by the functions $\cos^2$ and $\sin^2$ respectively.

This general result for the wavespeeds gives a dependence on the angle variable $\theta$.
We consider some general numerical experiments to assess the behaviour of the system.
What we find is two distinct types of behaviour.
We refine our wavespeed estimates in light of these particular cases.

\subsection{Small parameter behaviour}
For parameters close to those of the simple system, namely small convection ($\gamma \ll 1$) and near-equal diffusion ($\varepsilon_1 \approx \varepsilon_2$), we expect the behaviour to remain similar to that of the simple system.
Previously we have seen ODE-like and Turing-type behaviour in a convection-free system~\cite{Flach2007Turing-pattern-}.
In both cases, we see that the internal angle of the pattern (the group direction) is not in line with the wavefront.
As such, the oscillation is perceptible at the wavefront.
This suggests that an average value approximation of $\theta$ may be an appropriate estimate.

In this case we approximate $\sin^2 \theta$ with its average value: $1/2$.
This yields $\bar{c} \approx \bar{\gamma}/2 \pm 2\sqrt{1+\bar{\varepsilon}/2}$.
In the unscaled coordinate system we have  $c \approx \gamma/2 \pm\sqrt{2}\sqrt{\varepsilon_1+\varepsilon_2}$.

If we consider the limit as the amount of diffusion tends to zero, we see that the centreline of the propagation is at $\gamma/2$.
In this situation, the convection does not increase the separation of the left and right wavefronts; the spread of the propagation is solely due to diffusion.
In the original system the incident fronts are at $\textit{wavespeed} \approx \cfrac{p+q}{2} \pm\sqrt{2}\sqrt{\varepsilon_1+\varepsilon_2}$.

For dimensional reasons, it is appropriate to write $\varepsilon_i = \hat{\varepsilon_1}^2$.
Then the diffusive contribution to the wavespeed becomes $\sqrt{2}\sqrt{\hat{\varepsilon_1}^2+\hat{\varepsilon_2}^2}$, which is reminiscent of Pythagorean form.
More simply, the term is additive and the wavespeed is proportional to the diffusion of each species.

\subsection{Large parameter behaviour}
As the strength of convection $\bar{\gamma}$ increases, we see an alignment of the incident wave with the internal group direction.
This means that the value of $\theta$ is approximately constant at the onset of the wave.
We investigate the effect of this on the wavespeed estimate.

We consider the range of values possible for the right wavespeed with a fixed value for $\theta$.
Since $\sin^2$ has a range of $[0, 1]$ then $\bar{c}$ must range at least between $2$ and $\bar{\gamma} + 2\sqrt{1 + \bar{\varepsilon}}$.
This second value is not necessarily greater than the first, for example $\bar{\gamma} = 0$ and $\bar{\varepsilon} = -1/2$ gives $\sqrt{2} \approx 1.4$, which is less than $2$.
However, once $\bar{\gamma} > 2$ the second value is greater than the first, so we take this as an estimate for the transition to large parameter behaviour.
Thus the range of possible values for the right wavespeed is $\bar{c}_R \in [2, \bar{\gamma} + 2\sqrt{1 + \bar{\varepsilon}}]$.

On the left we have a similar result: $\bar{c}_L \in [-2, \bar{\gamma} - 2\sqrt{1 + \bar{\varepsilon}}]$ is the range once $\bar{\gamma}>2$.
We note that for $\sin^2=0$, $\theta = 0 \pmod{\pi}$ and for $\sin^2=1$, $\theta=\pi/2 \pmod{\pi}$.

We propose that the maximum spread of the instability is achieved.
That is, the wavespeeds are maximised within their given possible ranges.
For this, the left wavespeed must reach its lower limit, and the right its highest.
This gives $\bar{c}_L \approx -2$ and $\bar{c}_R \approx \bar{\gamma} + 2\sqrt{1+\bar{\varepsilon}}$.
Then the onset angle to the left is expected to be approximately $0 \pmod{\pi}$, and to the right of the right wave the angle should be roughly $\pi/2 \pmod{\pi}$.

\begin{figure}
\centerline{\includegraphics[width=4.2in]{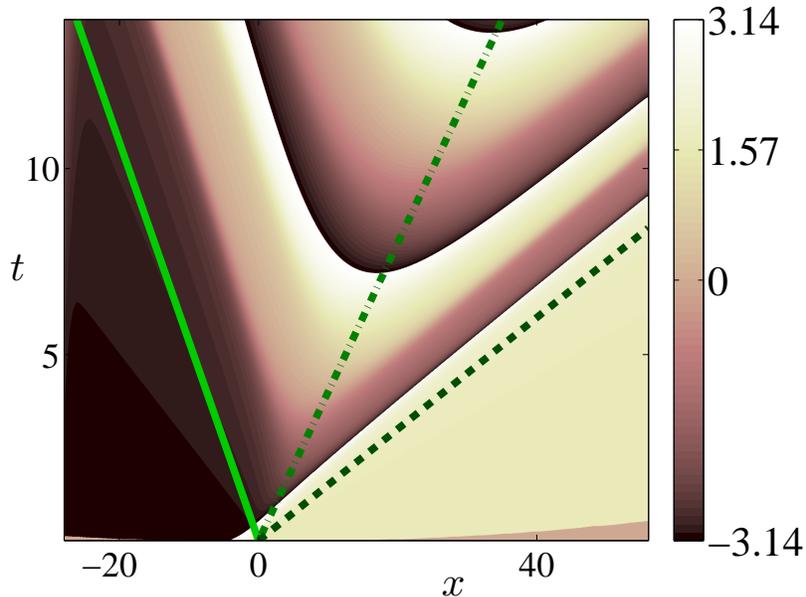}}
\caption{Large parameter behaviour of $\theta$.
The oscillations are aligned with the travelling waves.
To the left of the left travelling wave (the light, solid line), the angle is roughly $0$, and to the right of the right wave (the dark, dashed line) the angle is approximately $\pi/2$.
There is a smooth transition between the left and the right of the pattern.
The oscillations at the centre line ($\gamma/2$; medium, dot-dashed line) are roughly perpendicular to the centre line itself.
This is a numerical solution of~\eqref{eq:full} using NAG D03PCF, plotting species $u$ with $\gamma=5$, $\varepsilon_1=\varepsilon_2=1$.
The reactants are initially at steady state: $\left(u,v\right)=\left(0, 0 \right)$, with a small disturbance at $x=0$.
The boundaries are held at zero derivative: $u_x=0$, $v_x=0$.}\label{fig:polar_angle}
\end{figure}
We see this is the case in {\sc Figure}~\ref{fig:polar_angle}.
There is a smooth transition between the left and the right of the pattern.
The oscillations at the centre line (approximately $\gamma/2$) are roughly perpendicular to the centre line itself.

In the unscaled coordinate system, for strong convection, we have $c_L \approx -2\sqrt{\varepsilon_1}$ and $c_R \approx \gamma + 2\sqrt{\varepsilon_2}$.
In the original system, this is $\textit{left wavespeed} \approx p - 2\sqrt{\varepsilon_1}$ and $\textit{right wavespeed} = q + 2\sqrt{\varepsilon_2}$ for $q > p$.
In the other case $p > q$, there is a reversal and $\textit{left wavespeed} \approx q - 2\sqrt{\varepsilon_2}$ and $\textit{right wavespeed} = p + 2\sqrt{\varepsilon_1}$.

This analytical prediction was made from observing that the internal waves become aligned with the onset of the instability.
However, when we make the assumption that the spread of the propagation is maximised, it follows that the onset angle $\theta$ become constant.
Furthermore, the suggestion of $\bar{\gamma}=2$ as a transition value between behaviours occurs naturally from the algebra.

\section{Numerical experiments}
\begin{figure}
\subfigure[Left wavespeed. The darker surface is the low parameter estimate: $\gamma - \sqrt{2\varepsilon}$; the lighter surface is the high parameter estimate: $-2$. Both estimates hold well with a smooth transition between behaviours.]
{\includegraphics[width=3in]{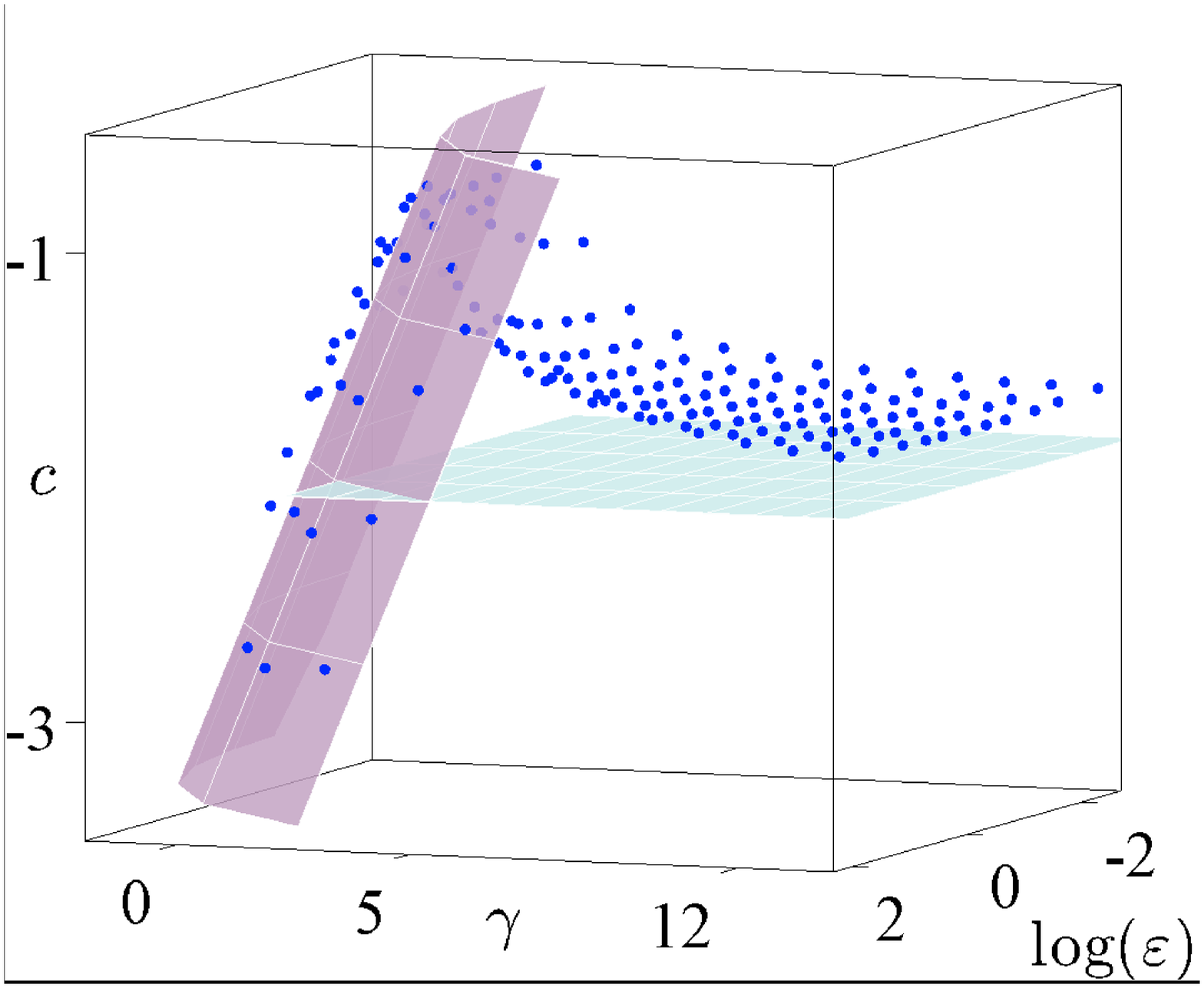}}
\subfigure[Right wavespeed. The darker surface is the low parameter estimate: $\gamma + \sqrt{2\varepsilon}$; the lighter surface is the high parameter estimate: $\gamma + 2\sqrt{\varepsilon}$. Both estimates are good and there is a neat transition between behaviours.]
{\includegraphics[width=2.82in]{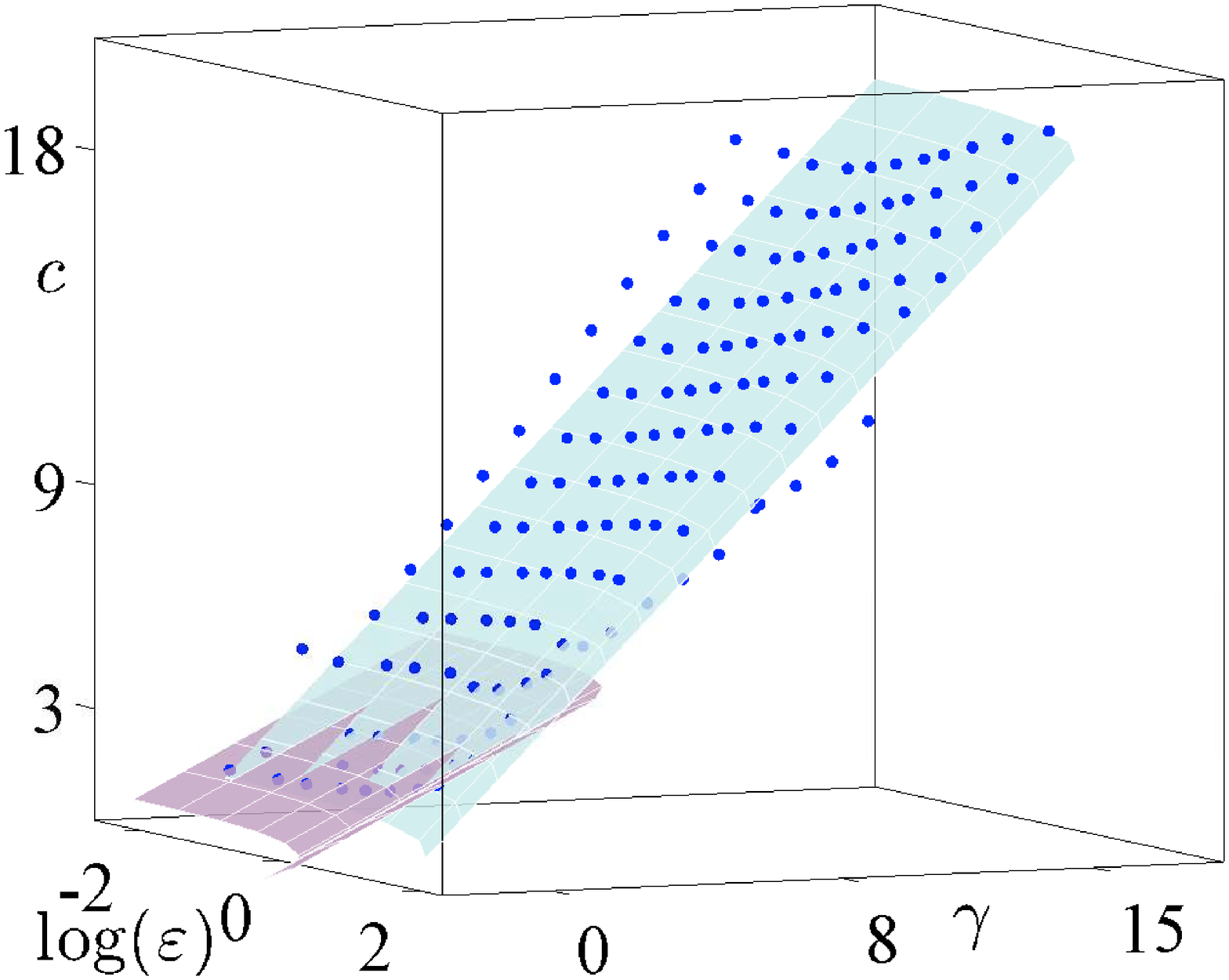}}
\caption{Numerical data compared to analytical estimates. The estimates correspond to the data over the full range. The low parameter estimate holds for longer for high $\varepsilon$ than for low.
The points are data from numerical solutions of~\eqref{eq:full} using NAG D03PCF, plotting species $u$.
The parameter $\varepsilon = \varepsilon_2/\varepsilon_1$; in this experiment $\varepsilon_1=1$ so $\varepsilon = \varepsilon_2$.
The reactants are initially at steady state: $\left(u,v\right)=\left(0, 0 \right)$, with a small disturbance at $x=0$.
The boundaries are held at zero derivative: $u_x=0$, $v_x=0$.}
\label{fig:ef05fig5}
\end{figure}
Having carried out an analysis together with some initial numerical investigation, we now conduct a more comprehensive experiment.
We vary our parameters over a wide range and measure the incident wavespeeds in each case.
The results confirm our analyses, as we can see in {\sc Figure}~\ref{fig:ef05fig5}.

For the left wavespeed the difference between the low and high parameter estimate is pronounced, so the transitional behaviour is noticable.
There is less difference between the estimates on the right, and the transition is close to the intersection of the estimates.
This seems quite remarkable, that the behavioural transition occurs where the estimates intersect.
As a result, there is no discernable area of transitional behaviour.

In both cases the low parameter behaviour endures longer for $\varepsilon$ high, compared to $\varepsilon$ low.
This suggests that the transition, caused primarily by the convection, $\gamma$, is held in check more strongly by the diffusion on $v$, where the convection is applied.

\section{Absolute versus convective instability}
There is a concept of stability in systems such as these, which considers whether the system at a particular spatial point, once destabilised, will return to the steady state.
If so, the system is classified as convectively unstable.
If the chemical species do not return to stasis then the system is deemed absolutely unstable.

In line with our analysis, we define absolute instability as $\sgn \left(\textit{left wavespeed} \right) \ne \sgn \left(\textit{right wavespeed} \right)$, with equivalence giving convective instability.
Combining this with our results, we see that for small parameters we require $\gamma^2 > 8 \left( \varepsilon_1+\varepsilon_2 \right)$ for convective instability in the unscaled system.
This inequality is unlikely to be satisfied whilst keeping within the small parameter regime.

However, if we consider the original system, the condition is $\left(p+q\right)^2 > 8 \left( \varepsilon_1+\varepsilon_2 \right)$, which is quite possible to satisfy while still remaining in the `small parameter' regime).
For example, $p=q=2$ gives $\gamma=0$, and then $\varepsilon_1= \varepsilon_2 = 1$ satisfies the condition.

For large parameters we always have absolute instability in the unscaled system.
Referring to the observed data ({\sc Figure}~\ref{fig:ef05fig5}), we see the left and right wavespeeds always have different signs.
Thus, we can classify this system as only exhibiting absolutely unstable behaviour.

However, if we consider the original system~\eqref{eq:original} we find that $q>p>2\sqrt{\varepsilon_1}$ or $p>q>2\sqrt{\varepsilon_2}$ will give convective instability rather than absolute instability.

We see that, in the original system, we can always shift the stability from absolute to convective with sufficient convection applied to the system.
In the flow coordinate system, the stability is essentially absolute.
However, when we convert to this system we make an arbitrary choice to remove the parameter p rather than, for example, q.
Then the distinction between absolute and convective instability becomes similarly arbitrary.

\section{An alternative perspective}
At the start of our analysis, in section~\ref{sec:pattern}, we chose a new coordinate system for the purpose of eliminating one of the convective terms.
This reduction in the number of parameters proved useful in the ensuing work.
However, the choice of coordinates, as we observed previously, was arbitrary.

Now we consider an alternative coordinate system, one centred on the convection.
From~\eqref{eq:original} we choose $x= \xi -\cfrac{p+q}{2} t$.
This gives the system as: \begin{eqnarray}
    u_t &=&  \varepsilon_1 u_{xx} + \hat{\gamma} u_x + f              \ , \nonumber\\
    v_t &=&  \varepsilon_2 v_{xx} - \hat{\gamma} v_x + g \ , \label{eq:flow}\end{eqnarray}
where $\hat{\gamma} = \cfrac{q-p}{2}$.
We can recalculate the apparent wavespeeds from this perspective.
For our `small parameter' regime we have $\hat{c} = \pm \sqrt{2}\sqrt{\varepsilon_1 + \varepsilon_2}$.
Here we see that we have aligned our system centrally with the flow.

In our large parameter case, we first consider $q>p$.
This ensures $\hat{\gamma} > 0$.
Now the left wavespeed is $\hat{c}_L = -\left(\hat{\gamma} + 2\sqrt{\varepsilon_1}\right)$, and the right wavespeed is $\hat{c}_R = \hat{\gamma} + 2\sqrt{\varepsilon_2}$.

In the other case $p<q$, we define $\tilde{\gamma} = \cfrac{p-q}{2} > 0$:\begin{eqnarray}
    u_t &=&  \varepsilon_1 u_{xx} - \tilde{\gamma} u_x + f              \ , \nonumber\\
    v_t &=&  \varepsilon_2 v_{xx} + \tilde{\gamma} v_x + g \ . \label{eq:flow2}\end{eqnarray}
Now the wavespeeds are  $\tilde{c}_L = -\left(\tilde{\gamma} + 2\sqrt{\varepsilon_2}\right)$, and the right wavespeed is $\tilde{c}_R = \tilde{\gamma} + 2\sqrt{\varepsilon_1}$.

Consider the right wavespeed.
In each case it is a function of only one of the diffusion parameters.
We examine the equation for the relevant species.
In both cases, we find a negative convection term for that species.
This suggests that negative convection drives the onset of instability.

For the left wavespeed the connection is with the positive convection.
However, in this direction $x$ is decreasing.
From this perspective, the signs of the convective terms change.
In this sense then, the connection is the same.

If we understand that the negative convection drives the instability wavefront, then we can consider only the associated species.
For example, looking for the right wavespeed in~\eqref{eq:flow2}, we see that $u$ is the relevant species.
We make the assumption that the other species $v$ is much smaller than $u$, and so disregard it.
Then we can carry out a travelling wave analysis for this single species.
The result is exactly the same as we found in our polar-form analysis, with a much simpler approach.
We can confirm this behaviour by observing that the initiation angles $0$ and $\pi/2$ correspond to $u$ and $v$ only, respectively.

We can go one step further in unravelling the puzzle.
If we simply separate the species -- perhaps we could call this a super-linearisation -- we can produce two wavespeeds, one from each species.
In our example~\eqref{eq:flow2}, for the right wavespeed we have the possibility of $\tilde{\gamma}+2\sqrt{\varepsilon_1}$ from the $u$ equation and  $-\tilde{\gamma}+2\sqrt{\varepsilon_2}$ from $v$.
Then we apply the `maximisation' principle suggested earlier.
For large $\tilde{\gamma}$, we can be sure that the wavespeed given from the $u$ equation will maximise the propagation.
Thus we can find the wavespeeds by applying two simple concepts to the problem.

The behaviour for small convection relies on some more complex way of combining the two species.
The wavespeed in this case is dependent on the strength of diffusion of both the species.
Therefore the polar-form approach given previously retains some merit.

\section{Discussion}
The addition of convection to the reaction-diffusion model has produced many patterns~\cite{Bamforth2000Modelling-flow-,Satnoianu2003Coexistence-of-}.
Advection, being the broader definition, has potential to produce even more.
Aspects of the pattern have been observed, such as the different behaviour to the right and the left of the free-forming pattern~\cite{Satnoianu1998Spatio-temporal}.

Previously analyses have considered a variety of models, most restricted in different ways.
No movement for one species~\cite{Satnoianu1998Interaction-bet}, equal diffusion and equal convection~\cite{Bamforth2000Modelling-flow-} and fixed diffusion-convection ratio~\cite{Satnoianu2000Non-Turing-stat} have all been considered.
Here, we have considered the most general diffusion and convection possible, as in more recent work~\cite{Satnoianu2003Coexistence-of-}.

Much previous work was devoted to the discussion of `flow and diffusion distributed structures'~\cite{Satnoianu2001Parameter-space}.
This is spatial pattern with no variation in time.
Our reduction of the system to two parameters allows for a clear and concise analysis of the system.
However, after our change to convective coordinates, the meaning of `time-constant' is lost.
In fact, in our reduced system, a time-constant solution is only possible for the trivial case, $\gamma=0$.

The standard analysis to date has been a linear existence analysis for pattern, involving dispersion relations.
This follows the Turing derivation, but is difficult for this more complex model.
In our analysis we have employed a different approach.
We have not directly considered the internal behaviour of the pattern, but instead looked for approximate descriptions of the propagation speed of the pattern.
This information has already been found for other versions of this system (less general diffusion or convection, more complex reactions)~\cite{Satnoianu1998Interaction-bet, Bamforth2000Modelling-flow-}.

To understand the understand a system well, it is best to isolate influences~\cite{Cross1993Pattern-formati}.
This is why we consider the pattern forming away from the spatial boundaries.
In this instance we also have symmetrical, parameterless reactions, giving the most elementary form possible.
Previous work has shown boundary effects mixed in with other behaviour~\cite{Satnoianu2003Coexistence-of-}, and that the boundary can control the existence of pattern~\cite{Bamforth2000Modelling-flow-}. 

Our analysis is the natural one given the system and its behaviour.
We have combined two simple techniques to reduce a complex problem to a more manageable one.
We consider only the (apparently) linear aspect.
Our analysis predicted the speed of the incident travelling waves over the full parameter range.
In combination with the numerical investigation, we demonstrate two different behaviours of the system.
There is a distinct transition between them, with the transition region indicated by the analysis.
Our work also confirms the validity of the Fisher-type analysis in this two-species system.
Furthermore, we give some insight to the different behaviour to right and left of the pattern.

Having successfully analysed this simple $\lambda$-$\omega$ system we must attempt a similar analysis on a more complex system.
The next step is to consider a more realistic set of reaction functions, such as the Schnakenberg ones.
In this case, we already know that first-order effects will come into play~\cite{Flach2006Limit-cycles-in}.
This should lead to a clearer understanding of the Turing phenomenon.

\section*{Acknowledgements}
We would like to acknowledge support from NIH grant number R01GM076692.
Any opinions, findings, conclusions or recommendations expressed in this paper are those of the authors and do not necessarily reflect the views of the NIH or the United States Government.
The authors are grateful to two anonymous reviewers for their helpful comments.

\bibliographystyle{prsty}

\end{document}